\documentclass[letterpaper, 10 pt, journal, twoside]{ieeetran}

\usepackage{graphicx}          

\usepackage{times,mathptmx}
\usepackage[T1]{fontenc}
\usepackage{color}
\usepackage{epsfig}
\usepackage[latin1]{inputenc}
\usepackage{rotating}
\usepackage{amsmath}
\usepackage{amssymb}
\usepackage{amsfonts}
\usepackage{epsf}
\usepackage{graphicx}
\usepackage{psfrag}
\usepackage{subfigure}
\usepackage{url}
\usepackage{tikz}
\usepackage{caption}
\usetikzlibrary{shapes,arrows}
\usepackage{pgfplots}
\pgfplotsset{plot coordinates/math parser=false}

\renewcommand{\baselinestretch}{0.89}

\newcommand{\R}{\mathbb{R}}

\title{Algorithms and Performance Analysis for Stochastic Wiener System Identification}

\author{\parbox{3 in}{\centering Bo Wahlberg*~\IEEEmembership{Fellow,~IEEE}
\thanks{This work was partially supported by the Swedish Research Council. }
\thanks{*Department of Automatic Control, School of Electrical Engineering and Computer Science, KTH Royal Institute of Technology, SE-100 44 Stockholm, Sweden.
 (e-mail: bo@kth.se)}}
  \hspace*{ 0.5 in}
  \parbox{3 in}{ Lennart Ljung**~\IEEEmembership{Fellow,~IEEE}
 \thanks{**Division of Automatic Control, Link\"oping University, SE-581 83 Link\"oping, Sweden}}
 }

\begin{document}

\maketitle
\thispagestyle{empty}
\pagestyle{empty}
\begin{abstract}   
We analyze the statistical performance of identification of stochastic dynamical systems with non-linear measurement sensors. This includes stochastic Wiener systems, with linear dynamics, process noise and measured by a non-linear sensor with additive measurement noise. There are many possible system identification methods for such systems, including the Maximum Likelihood (ML) method and the Prediction Error Method (PEM). The focus has mostly been on algorithms and implementation, and less is known about the statistical performance and the corresponding Cram\'er-Rao Lower Bound (CRLB) for identification of such non-linear systems.  We derive expressions for the CRLB and the asymptotic normalized covariance matrix for  certain Gaussian approximations of Wiener systems to show how a non-linear sensor affects the accuracy compared to a corresponding linear sensor. The key idea is to take second order statistics into account by using a common parametrization of the mean and the variance of the output process. This analysis also leads to a ML motivated identification method based on the conditional mean predictor and a Gaussian distribution approximation. The analysis is supported by numerical simulations.

\end{abstract}

\section{Introduction}
\IEEEPARstart{T}{here} has recently been a lot of progress in the development of algorithms for identification of nonlinear system, see \cite{7470636} for a recent overview. The focus has been on different deterministic model structures such as Best Linear Approximation, Volterra, Hammerstein and Wiener models and its generalizations. However, most of these methods assume rather restricted noise models, such as predictor models with additive white noise. Identification of Wiener systems, that is systems with a static nonlinearity at the output,  is a classical problem in system identification, see \cite{Ljung:99,Billing&Fakhouri:82,bai2003frequency,zhu2002estimation,enqvist2005linear,Zhao20121175,Bendat98}. It is a key component in the identification of block based non-linear systems, and we refer to the edited book \cite{Fouad&Bai:2010} for many recent contributions.  Identification of Wiener systems is an extensive topic and our reference list is by no means complete. However, most of the research and existing algorithms for Wiener system identification assume that the process noise can be neglected. The corresponding system identification optimization problem is then solved by minimizing the error between the measured and simulated outputs. As pointed out in \cite{hagenblad2008maximum}, process noise will then cause a biased estimate. The Maximum Likelihood (ML) method for stochastic Wiener models was first presented in \cite{914162}, and described in more  detail in \cite{hagenblad2008maximum}. Efficient numerical algorithms of the ML estimation problem  based on the EM algorithm and the particle filter have been presented in \cite{Schon201139,wills2013identification,wills2010wiener}. The statistical properties of the ML method for identification of stochastic Wiener systems are, however, less studied. The Cram\'er-Rao Lower Bound  for systems without process noise, where the linear part is a FIR model and the static nonlinearity is a polynomial, is derived in \cite{nordsjo1997cramer}. More recent results on stochastic Wiener system identification including benchmarks can be found in \cite{7828033,SCHOUKENS2017446,7039904,WAHLBERG2015620,7798727}.

The main contributions of the paper are:
\begin{itemize}
\item  Approximate expressions for the Fisher Information Matrix  and the asymptotic normalized covariance matrix  for identification of stochastic Wiener systems that give insight how a nonlinear sensor affects the accuracy of the identified model.
 \item An identification algorithm for stochastic Wiener systems based on the Conditional Mean Predictor and a Gaussian distribution approximation is derived and analyzed.
\end{itemize}

The structure of this paper is as follows:   Maximum Likelihood identification of nonlinear system is summarized  in Section \ref{sec:ML}. Formulas for the Fisher Information Matrix and  Cram\' er-Rao Lower Bound are presented in Section \ref{sec:FIM}. The special case of estimating the mean of a Gaussian process using a nonlinear sensor is analyzed in detail in Section~\ref{sec:SC}. Here the Fisher Information Matrix and the asymptotic normalized covariance matrix are derived for certain Gaussian approximations. This leads to asymptotic error variance results on how the nonlinearity affects the estimate. These results are then generalized to stochastic Wiener systems in Section \ref{sec:GSWS}. An example is outlined and  numerically evaluated in Section \ref{sec:IE}. The paper is concluded in Section \ref{sec:9}.

\section{The Maximum Likelihood Method}
\label{sec:ML}
This section summarizes some general results on system identification of stochastic non-linear systems and is mainly based on \cite{Schon201139}.
Consider the non-linear stochastic state-space model structure
\begin{align}
x_{t+1}&=f(x_t,u_t,v_t,\theta)\nonumber\\
y_{t}&=h(x_t,u_t,e_t,\theta),
\label{eq:sys1}
\end{align}
with the state-vector $x_t\in \R^{n}$, input-signal sequence $\{u_t \in \R\}$,  output-signal sequence $\{y_t \in \R \}$. The process noise $\{v_t\ \in \R\}$ and the measurement noise $\{e_t \in \R\}$ are assumed to be mutual independent i.i.d.~processes with probability density functions (pdf's) $p_v(\cdot)$ and $p_e(\cdot)$, respectively. The parameters to be estimated are the elements of the vector $\theta \in \R^m$. It is often convenient to represent the model (\ref{eq:sys1}) in the stochastic Markov form
\begin{align}
x_{t+1}&\sim p_\theta(x_{t+1}| x_t)\nonumber\\
y_{t}&\sim p_\theta (y_t|x_t ),
\end{align}
where the conditional pdf's describe the dynamics of (\ref{eq:sys1}).  Here we use the same symbol $p_\theta(\cdot ) $ for different pdf's and let its argument decide which function to use. Later, we will  use an extra sub-index to more clearly define a specific pdf.
We let $\theta=\theta_o$ denote the true data generating system.

The system identification problem is to estimate $\theta_o$ from $N$ measurements of the input-output response
\begin{equation}
U_{t}=[u_1,\ldots, u_t],\quad Y_{t}=[y_1,\ldots, y_t],\quad t=1\dots , N.
\end{equation}
In order to compute the likelihood  function, we apply the so-called measurement update
\begin{align}
p_\theta(y_t|Y_{t-1})&=\int  p_\theta (y_t|x_t) p_\theta (x_t|Y_{t-1})dx_t,
\label{eq:m update}
\end{align}
and the time update
\begin{align}
p_{\theta} (x_t|Y_{t})& =\frac{p_\theta(y_t|x_t)p_{\theta} (x_t|Y_{t-1})}{p_{\theta}(y_t|Y_{t-1})},\\
p_{\theta}(x_{t+1}|Y_t)&=\int  p_\theta(x_{t+1}|x_t) p_{\theta} (x_t|Y_{t})dx_t.\label{eq:t update}
\end{align}
The log-likelihood function $\log p_\theta(Y_n)$ then equals
\begin{align}
l_\theta(Y_N)&=\log p_\theta(y_1)+\sum_{t=2}^N\log p_\theta(y_t|Y_{t-1}).
\end{align}
The Maximum Likelihood (ML) estimate $\hat{\theta}_{ml}$ is obtained by maximizing the cost-function $l_\theta(Y_N)$ with respect to $\theta$.
\section{The Fisher Information Matrix  and the Cram\' er-Rao Lower Bound}
\label{sec:FIM}
Define the Fisher Information Matrix 
\begin{equation}
I_{\theta_o}(N)={\mathrm E}\left[\frac{\partial l_\theta(Y_N)}{\partial \theta}\frac{\partial l_\theta(Y_N)}{\partial \theta}^T\right]_{\theta=\theta_o}.
\label{FIMtrue}
\end{equation}
The covariance matrix ${\rm Cov}\{\hat{\theta}\}$ of  any unbiased estimator $\hat{\theta}$ of $\theta_o$ satisfies the Cram\' er-Rao bound
\begin{equation}
{\rm Cov}\{\hat{\theta} \}\geq [I_{\theta_o}(N)]^{-1}.
\end{equation}
In the case of scalar Gaussian distributed observations
\begin{equation}
y_t \sim \mathcal{N}\left(\mu_t(\theta),C_t(\theta)\right),\quad t=1,\ldots N,
\end{equation}
the parameter dependent part of the log-likelihood function equals
\begin{equation}
l_\theta(Y_N)=-\frac{1}{2}\left[\sum_{t=1}^N \frac{[y_t-\mu_t(\theta)]^2}{C_t(\theta)}+ \log C_t(\theta)\right].
\label{eq:negloglike}
\end{equation}
The corresponding  Fisher Information Matrix has the form, see \cite{kay1993fundamentals},
\begin{align}
I_{\theta_o}(N)=\sum_{t=1}^N \left[\frac{1}{C_t(\theta)}\frac{\partial \mu_t(\theta)}{\partial \theta}\frac{\partial \mu_t(\theta)}{\partial \theta}^T\right.\qquad\qquad\qquad \nonumber\\ \left. +\frac{1}{2} \frac{1}{C_t(\theta)^2}\frac{\partial C_t(\theta)}{\partial \theta} \frac{\partial C_t(\theta)}{\partial \theta}^T \right]_{\theta=\theta_o}.
\label{eq:FIM}
\end{align}
This result only holds for Gaussian distributed noise. A less  well known result is that the score covariance matrix (\ref{FIMtrue}) corresponding to the cost function (\ref{eq:negloglike}) but for general noise distribution equals
\begin{align}
J_{\theta_o}(N)= \sum_{t=1}^N \left[\frac{1}{C_t(\theta)}\frac{\partial \mu_t(\theta)}{\partial \theta}\frac{\partial \mu_t(\theta)}{\partial \theta}^T\right.\qquad\qquad\qquad \nonumber\\ \left. +\frac{\kappa(\theta)}{2}\frac{1}{C_t(\theta)^2} \frac{\partial C_t(\theta)}{\partial \theta} \frac{\partial C_t(\theta)}{\partial \theta}^T \right]_{\theta=\theta_o}.
\label{eq:FIMgeneral}
\end{align}
where
\begin{equation}
\kappa(\theta)=\frac{D_t(\theta)}{2C_t(\theta)^2},\: D_t(\theta)={\rm E}\{[(y_t-\mu_t(\theta))^2-C_t(\theta)]^2\}.
\label{eq:kap}
\end{equation}
The derivation of (\ref{eq:FIMgeneral}) is based on calculations of the gradient of (\ref{eq:negloglike}) as done in  Expression $(3C.6)$ in Appendix 3C in \cite{kay1993fundamentals}.
 Notice that $D_t(\theta)=2C_t(\theta)^2$ and thus $\kappa(\theta)=1$ for a Gaussian distribution, which gives back the result (\ref{eq:FIM}).  The  kurtosis of a stochastic process equals $2\kappa(\theta)+1$, and is a standard measure of infrequent extreme deviations from the mean of the process. For example, the kurtosis for a standard chi-squared distributed variable with one degree of freedom is $15$ compared to $3$ for the Gaussian case.
 The motivation for using (\ref{eq:negloglike}) for a non-Gaussian distribution is that it can be viewed as an extension of the standard Prediction Error Method (PEM) by also matching the second order statistics.

We will mainly be interested in the asymptotic (large $N$) performance of the identification methods, which is measured by the asymptotic normalized covariance matrix
\begin{equation}
{\rm AsCov}\{\hat{\theta}\}=\lim_{N\to \infty} {\rm Cov}\{\sqrt{N}[\hat{\theta}-\theta_o]\},
\end{equation}
and by the asymptotic Fisher Information Matrix (FIM) and corresponding asymptotic Cram\' er-Rao Lower Bound (CLRB)
\begin{equation}
{\rm FIM}(\theta_o)=\lim_{N\to \infty} \frac{1}{N}I_{\theta_o}(N), \quad {\rm CRLB}(\theta_o)=[{\rm FIM}(\theta_o)]^{-1}.
\label{eq:CRB}
\end{equation}
The ML method is under certain regularity conditions asymptotically efficient in the sense that it achieves the asymptotic CLRB, \cite{Ljung:99,kay1993fundamentals},
\begin{equation}
{\rm AsCov}\{\hat{\theta}_{ml}\}={\rm CRLB}(\theta_o).
\end{equation}
The cost-function (\ref{eq:negloglike}) makes sense even for a non-Gaussian distribution. The $\log$-part can be seen as a regularization term to penalize a too large variance estimate.
The asymptotic normalized covariance  matrix of the estimate obtained by minimizing (\ref{eq:negloglike}) for a general noise distribution equals
\begin{equation}
{\rm AsCov}\{\hat{\theta} \}=\lim_{N\to \infty} {N}[I_{\theta_o}(N)]^{-1}J_{\theta_o}(N)[I_{\theta_o}(N)]^{-1}.
\label{eq:asvarmod}
\end{equation}
It is only related to  the CRLB under the Gaussian assumption for which $J_{\theta_o}(N)=I_{\theta_o}(N)$.
Expression (\ref{eq:asvarmod}) follows, c.f.~Chapter 9.2 in \cite{Ljung:99}, from analyzing
 $$
 0=\frac{\partial l_\theta(Y_N)}{\partial \theta}|_{\theta=\hat{\theta}}
 \approx \frac{\partial l_\theta(Y_N)}{\partial \theta}|_{\theta=\theta_o}+\frac{\partial^2 l_\theta(Y_N)}{\partial \theta^2}|_{\theta=\theta_o}[\hat{\theta}-\theta_o].
 $$

{\bf To conclude:} The performance of identification methods for stochastic non-linear systems can in principle be evaluated using the results described in this section, for example by the CRLB and the asymptotic normalized covariance matrix. However, one has typically to resort to numerical calculations to determine these expressions. In particular, it is difficult to obtain insights in how a specific non-linearity will affect the identification accuracy. The objective of the paper is to give a more transparent results for the special case of identification of stochastic Wiener systems.

\section{Stochastic Wiener Systems}
\label{sec:SWS}
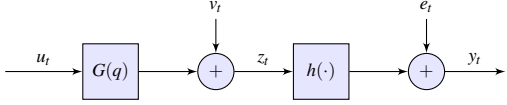
\begin{figure}[htb]
\centering
\pgfplotsset{width=1\columnwidth,height=0.5\columnwidth,compat=newest,plot coordinates/math parser=false}\tikzstyle{block} = [draw, rectangle,
    minimum height=3em, minimum width=3em,align=center,node distance=2cm,fill = blue!10]
\tikzstyle{sum} = [draw, circle, node distance=2cm,fill = blue!10]
\tikzstyle{input} = [coordinate]
\tikzstyle{output} = [coordinate]
\tikzstyle{corner} = [coordinate]
\tikzstyle{pinstyle} = [pin edge={to-,thin,black}]

\begin{tikzpicture}[auto, node distance=1cm,>=latex',transform shape,scale=0.7]

    \node [input] (input_u) {};
    \node [block,right of=input_u] (G) {$G(q)$};
    \node [sum, right of=G] (sum1) {$+$};
    \node [input,above of=sum1] (input_w) {};
    \node [block,right of=sum1] (f) {$h(\cdot)$};
    \node [sum, right of=f] (sum2) {$+$};
    \node [input, above of=sum2] (input_e) {};
    \node [output,right of=sum2,node distance=1.5cm] (output_y) {};
        
    \draw [draw,->] (input_u) -- node[above]{$u_t$} (G);
    \draw [draw,->] (G) -- node[above]{} (sum1);
    \draw [draw,->] (input_w) -- node[above,at start]{$v_t$} (sum1);
    \draw [draw,->] (sum1) -- node[above]{$z_t$} (f);
    \draw [draw,->] (f) -- (sum2);
    \draw [draw,->] (input_e) -- node[above,at start]{$e_t$} (sum2);
    \draw [draw,->] (sum2) -- node[above]{$y_t$} (output_y);
    
\end{tikzpicture}
\caption{Stochastic Wiener system.}\label{fig:1}
\end{figure}

Consider a stable scalar discrete time stochastic Wiener dynamic model structure illustrated in Figure~\ref{fig:1},
\begin{align}
z_t&=G(q,\theta)u_t+v_t,\nonumber\\
y_t&=h(z_t,\theta)+e_t,
 \label{eq:2}
\end{align}
with transfer function $G(q,\theta)=\sum_{k=0}^\infty g_k q^{-k},$
impulse response sequence $\{g_k \in \R\}$, ($q$ is the shift operator),  white zero mean process noise
$\{v_t \in \R\}$ with pdf $p_{v,\theta}(\cdot)$,  and additive zero mean white measurement
noise $\{e_t \in \R\}$ with pdf $p_{e,\theta}(\cdot)$.

It is possible to consider colored process noise $v_t$  by using a predictor form model  as described in Section 3.3 in \cite{hagenblad2008maximum}.
For a stochastic Wiener system the measurable output $y_t$ is a non-linear function $h(\cdot)$ of the output $z_t$ of a stochastic linear system. This can correspond to a non-linear sensor.

The filtering equations (\ref{eq:m update}) for stochastic Wiener models (\ref{eq:2}) is simplified since there is no correlations in time due to the white noise assumptions, i.e.~$p_\theta(y_t|Y_{t-1})=p_\theta(y_t)$ with
\begin{align}
p_\theta(y_t)&= \int  p_{e,\theta}(y_t-h(G(q,\theta)u_t+v,\theta)) p_{v,\theta}(v)dv\nonumber\\
& ={\mathrm E}_v\{ p_{e,\theta}(y_t-h(G(q,\theta)u_t+v,\theta))\}.
\label{eq;24a}
\end{align}
The interpretation is to marginalize (calculate the average)  of the pdf $p_{e,\theta}(y_t-h(G(q,\theta)u_t+v_t,\theta))$ with respect to the process noise $v_t$.
If the function $h(x)$ is invertible, $h^{-1}(h(x))=x$, an equivalent likelihood expression is
\begin{equation}
p_\theta(y_t)=\int \frac{p_{v,\theta}(h^{-1}(y_t-e)-G(q,\theta)u_t)}{|h'(h^{-1}(y_t-e))|} p_{e,\theta}(e)de,
\label{eq:24}
\end{equation}
where prime denotes the derivative. If $h(x)$ is increasing, the absolute value in the denominator can be removed. The integral (\ref{eq:24}) may be easier to compute numerically than (\ref{eq;24a}). Here the $\theta$ dependence of $h(\cdot,\theta)$ has been omitted due to notation constraints. This expression is obtained by change of integration variable $e=y_t-h(G(q,\theta)u_t+v)$. Notice that if $z=h(v)$, where the pdf of $v$ is $p_v(x)$ and $h$ is an strictly monotonic differentiable function,  then the pdf of $z$ equals
\begin{equation}
p_z(x)=\frac{p_v(h^{-1}(x))}{\left|\frac{d}{dx} h^{-1}(x))\right|}= \frac{p_v(h^{-1}(x))}{|h'(h^{-1}(x))|}.
\end{equation}
The main challenge of ML identification is  how to efficiently calculate the pdf integrals. Notice that one has to calculate one integral per measurement (in total $N$) just to evaluate  the log-likelihood cost function
$l_\theta(Y_N)$ at a certain value of $\theta$. It is also, in principle, possible to  numerically calculate the corresponding FIM and CRLB. An efficient  way to numerically calculate the expected value of a function of a Gaussian variable  is to use the Gauss-Hermite  Quadrature
\begin{equation}
\int_{-\infty}^\infty g(x)e^{-x^2}dx\approx \sum_{i=1}^nw_ig(x_i), \label{eq:ghq}
\end{equation}
where $x_i$ are the roots of the Hermite polynomial $H_n(x)$ of degree $n$ and the weights are $w_i={2^{n-1}n!}/[{n^2[H_{n-1}(x_i)]^2}]$.
The  Gauss-Hermite approximation is exact for polynomials $g(x)$ up to order $2n-1$. See \cite{1499276} for a survey on its use in non-linear filtering and the connection to the unscented transform.
It is more difficult to obtain insight, {\em e.g.}, in how the properties of $h(\cdot)$ affect the accuracy of the ML estimate of linear dynamics.

As mentioned in the introduction, PEM is an alternative to ML. For the stochastic Wiener system (\ref{eq:2}) the  conditional mean predictor equals
\begin{align}
\label{eq:LL1}
\hat{y}_{t|t-1}(\theta)&={\mathrm E}_v\{h (G(q,\theta)u_t+v),\theta)\},
\end{align}
which in general can be computed by integration over $v$.
In some cases it is possible to analytically calculate
$\hat{y}_{t|t-1}(\theta)$ and its variance.  This is, for example, the case when  $v_t$ is
Gaussian distributed and $h(x)$ is a polynomial.

\section{Special Case}
\label{sec:SC}
We will motivate algorithms and analysis for ML and PEM identification of stochastic Wiener systems by first studying the special case of estimating a scalar constant  $m_o$  from the measurement model
\begin{equation}
y_t=h(m+v_t)+e_t.
\label{eq:sc}
\end{equation}
This corresponds to a stochastic Wiener model with $G(q,\theta)=m$, $\theta=m$ and $u_t=1$. We assume $\{e_t\}$ and  $\{v_t\}$ to be independent zero mean Gaussian distributed stochastic processes with  given  variances $\sigma^2_e$ and $\sigma^2_v$, respectively.  We also assume that the sensor function $h(\cdot)$ is given and known.
\vskip\baselineskip
We now want to answer the following two questions:
\vskip\baselineskip
\begin{enumerate}
\item How does the quality of the ML estimate $\hat{m}$ depends on the possible non-linear function $h(\cdot)$ and the noise variances $\sigma_v^2$ and $\sigma_e^2$? More precisely,  how does the  asymptotic Fisher Information Matrix, the corresponding Cram\'er-Rao Lower Bound and the asymptotic normalized covariance matrix depend on $h(\cdot)$, $\sigma_v^2$ and $\sigma_e^2$?

\vskip 0.5\baselineskip

\item How should the sensor  $h(\cdot)$ be designed to suppress noise and at the same time amplify information about the unknown parameter?

\end{enumerate}

\subsection{FIM and CRLB Expressions}
\label{sub:A}
It is in principle possible to numerically calculate the FIM and CRLB using the formulas  given the previous sections. The corresponding results are, however, quite involved and we will instead derive some approximative FIM and CRLB expressions  for the model (\ref{eq:sc}).

To start, what can be learned from the linear case $h(x)=Kx$? The corresponding model is
$
y_t=Km+e_t+Kv_t,
$
with the ML estimate
\begin{equation}
\hat{m}_{ML}=\frac{1}{KN}\sum_{t=1}^Ny_t.
\end{equation}
The asymptotic variance of the scaled error $\sqrt{N}[\hat{m}_{ML}-m_o]$ is equal to 
\begin{align}
{\rm CRLB}(m_o)&=\frac{\sigma_e^2+h'(m_o)^2\sigma_v^2}{h'(m_o)^2} \quad \Rightarrow\label{eq:CRBlin}\\
 {\rm FIM}(m_o)& =\frac{h'(m_o)^2}{\sigma_e^2+h'(m_o)^2\sigma_v^2},\quad h'(m_o)=K.
\end{align}
Another special case where it is possible to analytically calculate the asymptotic CRLB is when $h(\cdot)$ is a general differentiable function but there is no process noise, {\em i.e.}~ $v_t=0$. For this case, \cite{kay1993fundamentals},
\begin{equation}
{\rm CRLB}(m_o)=\frac{\sigma_e^2}{h'(m_o)^2}.
\end{equation}
The proof is based on a Taylor series approximation. If instead $e_t=0$ and the function $h(\cdot)$ is invertible we have the relation
\begin{equation}
h^{-1}(y_t)=m+v_t,\Rightarrow\; {\rm CRLB}(m_o) =\sigma_v^2.
\end{equation}
By comparing these two non-linear special cases  with  CRLB  for the linear sensor (\ref{eq:CRBlin}),
\begin{equation}
\frac{\sigma_e^2}{h'(m_o)^2}+\sigma_v^2,
\label{eq:crlblin}
\end{equation}
we note that for low or high values of $\sigma_v^2$ relative to $\sigma_e^2$, we expect the CRLB for a  nonlinear sensor to be close to the CRLB  for the linear case (\ref{eq:crlblin}). What happens in between these two extremes is an open problem to be addressed.

\subsection{First Order Approximations}
Applying Gauss Approximation Formula
\begin{equation}
h(m+v_t)\approx h(m)+h'(m)v_t
\end{equation}
to (\ref{eq:sc}) gives the "first-order" Gaussian  model
\begin{equation}
y_t=h(m)+h'(m)v_t+e_t.
\label{eq:gaf}
\end{equation}
Notice that this is an approximative model and will not lead to a direct approximation of the CRLB for the general non-linear case. This approximation is usually only valid in a neighbourhood where $h(\cdot)$ is approximately linear.
Estimation of $m$ using the model (\ref{eq:gaf}) is a Gaussian identification problem with jointly parameterized  mean and  variance functions
\begin{equation}
\mu(m)=h(m),\quad C(m)=\sigma_e^2+h'(m)^2\sigma_v^2.
\label{eq:muc}
\end{equation}
Notice that the noise variances are assumed to be known, since we otherwise may have identifiability problems. We can now directly use the FIM expression (\ref{eq:FIM}), which for this case simplifies to
\begin{align}
{\rm FIM}(m_o)&=\frac{\mu'(m_o)^2}{C(m_o)}+\frac{1}{2}\frac{C'(m_o)^2}{C(m_o)^2}.
\label{eq:FmuC}
\end{align}
Evaluation of (\ref{eq:FmuC}) using (\ref{eq:muc}) gives the   FIM expression:
\vskip\baselineskip
{\bf Result 1:} The Fisher Information Matrix for the first order  Gauss approximation model (\ref{eq:gaf}) equals
\begin{align}
{\rm FIM}_1(m_o)&= \frac{h'(m_o)^2}{\sigma_e^2+h'(m_o)^2\sigma_v^2 }+2\left[\frac{\sigma_v^2h'(m_o)h''(m_o)}{\sigma_e^2+h'(m_o)^2\sigma_v^2}\right]^2.
\label{eq:FIMfo1}
\end{align}
This result gives several new insights. The first term of the FIM expression (\ref{eq:FIMfo1})
equals the FIM (\ref{eq:CRBlin}) for the linear case. Notice that this now holds for a non-linear model (\ref{eq:gaf}). The reason is that the Gauss approximation formula linearizes the non-linear noise contribution in an appropriate way. The second term of (\ref{eq:FIMfo1}) shows that the uncertainty of the ML estimate is further reduced by utilizing the $m$ dependence of the variance. The improvement in information is proportional to $h''(m_o)^2$, which makes sense since the variance depends on $h'(m)$.

The FIM expression (\ref{eq:FIMfo1}) shows that a non-linear sensor, with $h''(m_o)\neq 0$ and gain $h'(m_o)$, can give a more accurate estimate of $m_o$ than a linear sensor with the same gain $K=h'(m_o)$. The improvement is, however, in general moderate since $\sigma_v^2$ needs to be small for this approximation to hold.

The corresponding result for the special case when the measurement noise is small and the function $h(x)$ is invertible is slightly more involved. Applying Gauss Approximation Formula to
\begin{equation}
h^{-1}(y_t)=h^{-1}(h(m+v_t)+e_t)\approx m+v_t+\frac{1}{h'(m+v_t)}e_t.
\label{eq:hin}
\end{equation}
leads to the Gaussian  model
\begin{equation}
 h^{-1}(y_t)=  m+v_t+\frac{1}{h'(m)}e_t.
 \label{eq:gaf2}
\end{equation}
Here we have  also approximated the factor $h'(m+v_t)$ by $h'(m)$.  The stochastic process (\ref{eq:gaf2}) has mean and variance
\begin{equation}
\mu(m)=m,\quad C(m)=\sigma_v^2+\frac{\sigma_e^2}{h'(m)^2},
\label{eq:muc3}
\end{equation}
and using (\ref{eq:FmuC}) gives the FIM expression:
\vskip\baselineskip
{\bf Result 2:} The Fisher Information Matrix for the first order  Gauss approximation model (\ref{eq:gaf2}) equals
\begin{align}
{\rm FIM}_2(m_o)&=\frac{1}{\sigma_v^2+\sigma_e^2/h'(m_o)^2}+2\left[\frac{\sigma_e^2 h''(m)}{h'(m)( \sigma_e^2+h'(m_o)^2\sigma_v^2)}\right]^2.
\label{eq:FmuC3}
\end{align}
This results is quite similar to Result~1, (\ref{eq:FIMfo1}). The only difference is that $\sigma_v^2h'(m_o)$ is replaced by $\sigma_e^2/h'(m_o)$. The two FIM expressions (\ref{eq:FIMfo1}) and (\ref{eq:FmuC3}) are equal if $$\sigma_e^2=h'(m_o)^2\sigma_v^2,$$ which makes sense from a noise contribution point of view.


\subsection{Second Order Approximation}
A more accurate approximation of (\ref{eq:sc}) is the model
\begin{equation}
y_t=h(m)+h'(m)v_t+\frac{h''(m)}{2} v^2_t+e_t.
\label{eq:x}
\end{equation}
This leads to a rather complicated ML problem due to the chi-squared distributed noise $v_t^2$.  The model (\ref {eq:x}) is exact for a quadratic sensor function $h(\cdot)$ and otherwise an approximation. This approach is related to the unscented transform as discussed in e.g.~\cite{4518435}. Rewrite the model (\ref{eq:x}) as follows to obtain a zero mean noise contribution
\begin{equation}
y_t=h(m)+\frac{h''(m)\sigma_v^2}{2}+ h'(m)v_t+ \frac{h''(m)}{2} [v^2_t- \sigma_v^2 ]+e_t.
\label{eq:soa}
\end{equation}
The variance of $v^2_t- \sigma_v^2 $ for a Gaussian process equals $2\sigma_v^4$. Hence the approximation (\ref{eq:soa}) has  mean and variance
\begin{align}
\mu(m)&=h(m)+\frac{h''(m)\sigma_v^2}{2},\label{eq:mean2}\\
C(m) &= \sigma_e^2+h'(m)^2\sigma_v^2+h''(m)^2\frac{\sigma_v^4}{2}.
\end{align}
Applying  (\ref{eq:FmuC}) to  this model structure gives:
\vskip\baselineskip
{\bf Result 3:} The Fisher Information Matrix for the (second order) Gaussian approximation of the model (\ref{eq:soa}) equals
\begin{align}
{\rm FIM}_3(m_o)=& \frac{[h'(m_o)+h'''(m_o)\sigma_v^2/2]^2}{\sigma_e^2+h'(m)^2\sigma_v^2+h''(m)^2\sigma_v^4/2 }\nonumber \\ & +2\left[\frac{h''(m_o)\sigma_v^2(h'(m_o)+h'''(m_o)\sigma_v^2/2)}{\sigma_e^2+h'(m)^2\sigma_v^2+h''(m)^2\sigma_v^4/2}\right]^2.
\label{eq:FIMfo2}
\end{align}
The main difference compared to the first order FIM (\ref{eq:FIMfo1}) is the influence of the third order derivative $h'''(m_o)$. Possible improvement due to increased information depends on the the size and sign of this term.

As noted earlier the model (\ref{eq:soa}) is not Gaussian  and Result~3 has to be modified as described by (\ref{eq:asvarmod}) in order to obtain the asymptotic normalized covariance matrix of the estimate obtained by minimizing (\ref{eq:negloglike}). We use (\ref{eq:FIMgeneral}) to calculate
\begin{align}
J_3(m_o)&= \frac{[h'(m_o)+h'''(m_o)\sigma_v^2/2]^2}{\sigma_e^2+h'(m)^2\sigma_v^2+h''(m)^2\sigma_v^4/2 }\nonumber \\ & +2\kappa(m_o)\left[\frac{h''(m_o)\sigma_v^2(h'(m_o)+h'''(m_o)\sigma_v^2/2)}{\sigma_e^2+h'(m)^2\sigma_v^2+h''(m)^2\sigma_v^4/2}\right]^2.
\label{eq:j3}
\end{align}
where $\kappa(m_o)$ is defined by (\ref{eq:kap}).

\vskip0.5\baselineskip

{\bf Result 4:} The asymptotic normalized covariance  of the estimate obtained by minimizing (\ref{eq:negloglike}) for the model (\ref{eq:soa}) equals
\begin{equation}
{\rm AsCov}\{\hat{m}\}=\gamma(m_o) \frac{1}{{\rm FIM}_3(m_o)}, \quad \gamma(m_o)=\frac{J_3(m_o)}{{\rm FIM}_3(m_o)}.
\label{eq:asvarmod2}
\end{equation}
The scaling $\gamma$ is related $\kappa$ in (\ref{eq:j3}) and gives a measure of how worse the accuracy is compared to the Gaussian CRLB based on (\ref{eq:FIMfo2}).

\subsection{Conditional Mean Predictor Model}
The model (\ref{eq:soa}) has several interesting interpretations. For a cubic  sensor function  $h(m)$ it gives the conditional mean predictor of $y_t=h(m+v_t)+e_t$. However, the corresponding prediction error
is {\em not} Gaussian distributed. The PEM  framework developed in \cite{7039904} uses
\begin{equation}
\hat{y}(m)={\mathrm E}_v \{ h(m+v)\},\quad \sigma_\epsilon^2(m)={\mathrm E}\{[y_t-\hat{y}(m)]^2\},
\end{equation}
and analyze the estimate obtained by minimizing a variance weighted PEM cost-function. A more accurate approach to estimate $m$ is to use the Gaussian ML cost-function (\ref{eq:negloglike})  and the model
\begin{equation}
y_t=\hat{y}(m)+\epsilon_t(m),\; 
\end{equation}
where $\epsilon_t(m)$ is zero mean with variance $ \sigma_\epsilon^2(m)$. Thus $\mu(m)=\hat{y}(m)$ and $C(m) =\sigma_\epsilon^2(m)$.
This leads in general to a more accurate estimate than the weighted PEM since the parameter dependence of the variance is taken into account. The corresponding asymptotic covariance matrix is given by (\ref{eq:asvarmod2}).

\section{Stochastic Wiener Model}
\label{sec:GSWS}
The results in the preceding section can  be generalized to the stochastic Wiener model (\ref{eq:2}) with a given non-linear sensor $h(\cdot)$.   The simplest case would be to use  the Gauss approximation model
\begin{align}
z_t(\theta)&=G(q,\theta)u_t,\nonumber\\
y_t&=h(z_t(\theta))+h'(z_t(\theta))v_t+e_t,
 \label{eq:2app}
\end{align}
for which it is possible to directly apply the ML method (\ref{eq:negloglike}) using
\begin{align}
\mu_t(\theta)= h(z_t(\theta)),\quad C_t(\theta)= \sigma_e^2+h'(z_t(\theta))^2\sigma_v^2.
\end{align}
The corresponding asymptotic Fisher Information Matrix is again obtained by taking the average of (\ref{eq:FIM}). The derivative of $h(z_t(\theta))$ will play the same role as for the simple case even if the formulas will be more involved. Also in this case a non-linear sensor together with the modelling of the variance can improve the accuracy of the $\theta$ estimate compared to using a linear sensor.
The key question is still how valid the approximative  model is for the intended use.

A more accurate description of (\ref{eq:2}) is the conditional mean predictor model
\begin{align}
\label{eq:LL1app2}
\hat{y}_{t}(\theta)&={\mathrm E}_v\{h (G(q,\theta)u_t+v),\theta)\},\nonumber\\
y_t&=\hat{y}_{t}(\theta)+\epsilon_t(\theta),\quad \sigma_{\epsilon,t}^2(\theta)={\mathrm E}\{[y_t-\hat{y}_t(\theta)]^2\}.
\end{align}
The  model parameter $\theta$ is estimated by maximizing the  Gaussian log-likelihood (\ref{eq:negloglike}) using
\begin{align}
\mu_t(\theta)=\hat{y}_{t}(\theta) ,\quad C_t(\theta)= \sigma_{\epsilon,t}^2(\theta).
\end{align}
This approach has recently been studied in \cite{Abdalmoaty2017}. It is very efficient from an implementation point of view compared to the true ML method.

\section{Examples}
\label{sec:IE}
We will now in more detail study the problem how to estimate the mean $m$ from observations of
$
y_t=h(m+v_t)+e_t
$.
Notice that this is a special case of the stochastic ML system
\begin{align}
z_t&=\theta u_t+v_t\nonumber\\
y_t&=h(z_t)+e_t,
\end{align}
with $m=\theta$ and $u_t=1$. Furthermore, we will assume that the gain is normalized to $h'(m_o)=1$, with $m_o=1$ and that the noises have equal power, $\sigma_e=\sigma_v$. This is the special case when we according to our theory could benefit from a non-linear sensor.
For a quadratic sensor $h(x)=x^2/2$, we obtain the mean and variance
\begin{align}
\mu(m)&=\frac{m^2}{2}+\frac{\sigma_v^2}{2},\quad C(m)=m^2\sigma_v^2+\frac{\sigma^4_v}{2}+\sigma_e^2.
\end{align}
For a cubic sensor $h(x)=x^3/3 $, we have
\begin{align}
\mu(m)&=\frac{m^3}{3}+ m\sigma_v^2,\quad C(m)=m^4\sigma_v^2+4m^2\sigma^4_v+\frac{5}{3}\sigma_v^6+\sigma_e^2.
\label{eq:cub}
\end{align}
We can now calculate the FIM (\ref{eq:FIMfo2}) and the asymptotic normalized variance  (\ref{eq:asvarmod2}) of the model parameter $m$
to evaluate identification performance.
For example, the quadratic sensor with $m_o=1$ and $\sigma_e=\sigma_v$ gives
\begin{align}
{\rm FIM}_3 &=\frac{1}{2\sigma_v^2+0.5\sigma^4_v}+\frac{1}{2}\frac{4\sigma_v^2}{(2\sigma_v^2+0.5\sigma^4_v)^2},\nonumber\\
J_3 &=\frac{1}{2\sigma_v^2+0.5\sigma^4_v}+\frac{\kappa}{2}\frac{4\sigma_v^2}{(2\sigma_v^2+0.5\sigma^4_v)^2}.
\end{align}
The scaling $\kappa$, for the quadratic sensor case, ranges from $1$ to $2.2$ when $\sigma_v$ tends from $0$ to $1$. A more accurate performance measure here than the FIM is the asymptotic normalized variance
\begin{equation}
{\rm AsCov}\{\hat{m}\}=\gamma \frac{1}{{\rm FIM}_3}, \quad \gamma=\frac{J_3}{{\rm FIM}_3}
\label{eq:ga}
\end{equation}
For the quadratic sensor $h(x)=x^2/2 $ the scaling $\gamma$ varies from $1$ to $1.5$ in our example.
This means that we then can expect worse identification accuracy (larger variance) for higher noise variances than predicted by ${\rm FIM}_3$.
For the cubic sensor $h(x)=x^3/3 $ the scaling factor $\gamma$ is in the order of $1$ to $8$, which shows that the true noise distribution here can be far from Gaussian.

Next we will compare the asymptotic performance results with finite data simulations.
We use $N=1000$ observations in the numerical study. The ML estimate is obtained by minimizing
\begin{equation}
l(\theta)=\sum_{t=1}^N -\log {\mathrm E}_{\bar{v}}\{\frac{{\scriptstyle  \sigma_v}}{\scriptstyle  \sqrt{2\pi}} e^{-\frac{1}{2\sigma_e^2}[y_t-h(\theta u_t+\sigma_v\bar{v}))]^2}\}
\label{eq:eML}
\end{equation}
 calculated using Gauss-Hermite Quadrature of order $10000$. The standard deviation of the $ML$ estimate is calculated from $250$ noise realizations. We assume that $m>0$ to avoid the obvious identifiability problem using $h(x)=x^2/2$.

\begin{table}[h]
\begin{center}{\small
 \begin{tabular}{||l|| c |c| c| c|c||}
 \hline
$\sigma_v^2$ & $0.1$ & $0.25$ &  $0.5$ &  $0.75$ & $1$ \\ [0.5ex]
 \hline
Linear & $0.0141 $ &  $0.0224 $  &  $0.0316 $ &   $0.0387 $ & $0.0447$\\
Quadratic  &  $0.0105 $  & $0.0179 $ &  $0.0280 $  &  $0.0371 $  &  $0.0461$ \\
ML2 &  $0.0132$ &  $0.0219 $ & $0.0314$ & $  0.0452$ &   $ 0.0478$\\
Cubic &  $ 0.0133 $  & $ 0.0239 $ &  $ 0.0396 $  &  $ 0.0532 $  &  $0.0660 $ \\
ML3 &  $0.0129$ &  $0.0214 $ & $0.0281$ & $0.0330$ &   $0.0449$\\
\hline
 \end{tabular}}
 \vskip 0.2cm
\caption{Linear  denotes the asymptotic normalized standard deviation for $h(x)=1$.
Quadratic denotes the asymptotic normalized standard deviation for $h(x)=x^2/2$.
ML2 denotes the standard deviation obtained by minimizing (\ref{eq:eML}) for for $h(x)=x^2/2$.
Cubic denotes the asymptotic normalized standard deviation for $h(x)=x^3/3$.
ML3 denotes the standard deviation obtained by minimizing (\ref{eq:eML})  for $h(x)=x^3/3$.}
\label{tab:1}
\end{center}
\end{table}
Table \ref{tab:1} shows that the asymptotic results are in good agreement with the simulations. They provide  reasonable  estimates of the accuracy of the true ML method. Notice however
that the asymptotic standard deviations corresponds to minimizing (\ref{eq:negloglike}) and not to the true ML function (\ref{eq:eML}) as in the simulations results ML2 and ML3. The proposed conditional mean predictor based identification method obtained by maximizing the Gaussian  log-likelihood (\ref{eq:negloglike}) gives results that are very close to the asymptotic normalized standard deviations Quadratic and and Cubic and are thus not reported in the table. The performance of conditional mean predictor based identification method is quite close to true ML method, while the computational efficiency is many magnitudes faster.
The identification accuracy is slightly worse for the cubic sensor and rather high noise levels. A reason is that the approximate model (\ref{eq:cub}) has problems to capture the non-linear stochastic behavior due to the cubic function.

\section{Conclusion}
\label{sec:9}
System identification with linear sensors and linear dynamics is a very advanced subject with powerful tools for performance analysis. We have taken a step towards understanding more about performance analysis for identification of stochastic non-linear systems.  In order to obtain transparent results to give guidelines for e.g.~experiment design, we have studied the special case of a known non-linear sensor and an unknown stochastic linear system with both process noise and measurement noise. This is a difficult problem, even when the objective is just to estimate a constant mean. We have derived first order and second order Gaussian approximations for the Fisher Information Matrix and the asymptotic normalized covariance matrix. We started by asking the two questions:
\begin{enumerate}
\item How does the quality of the model estimate  depends on the possible non-linear sensor function $h(\cdot)$ and the noise variances?

\item How should the sensor function  $h(\cdot)$ be designed to enhance information about the unknown parameter?

\end{enumerate}

Our main focus has  been on the first question. We have derived explicit FIM and variance expressions that shows how  the higher derivatives of the non-linear sensor function affect the quality of the estimate.
This is a first step to understand how general nonlinear functions affect the accuracy of parameter estimates.  We have showed how the scaling $\kappa$, which is directly related to the kurtosis, may amplify or reduce uncertainty. From a practical point of view one should design sensor systems that results in a low kurtosis. For example, the case when $v_t$ is uniformly distributed  gives a smaller kurtosis than for the Gaussian  case. It would be of interest to study more complicated systems including stochastic neural networks based models. Most results in this area are deterministic and little is known for example how stochastic input noise affects the performance.

We have proposed an identification method based on the conditional mean predictor and the corresponding prediction error variance combined with the Gaussian ML cost-function. This  results in a very computational efficient algorithm with good performance for identification of stochastic Wiener systems.

\vskip0.5\baselineskip

\bibliography{refs,referenslista}

\end{document}